\documentclass[10pt]{article}
\usepackage{ijcai20}

\usepackage{times}
\usepackage{amsmath,amssymb,amsfonts,amstext,amsthm,mathrsfs}
\usepackage{enumitem}
\usepackage{thm-restate}
\usepackage{algorithmic,algorithm}
\usepackage{wrapfig}
\usepackage[algo2e]{algorithm2e}
\usepackage{cleveref}
\usepackage{booktabs} 
\usepackage{subcaption}
\usepackage{graphicx} 
\usepackage[font=small]{caption}

\newtheorem{definition}{Definition}

\newtheorem{assum}{Assumption}

\newcommand{\dist}{\mathrm{dist}}
\newcommand{\prox}[1]{\mathrm{prox}_{#1}}

\newcommand{\dom}{\mathop{\mathrm{dom}}}
\newcommand{\zero}{\mathbf{0}}
\newcommand{\crit}{\mathop{\mathrm{crit}}}
\newcommand{\argmin}{\mathop{\mathrm{argmin}}}
\newcommand{\KL}{K{\L}~}

\newcommand{\inner}[2]{\langle #1, #2 \rangle}

\title{Proximal Gradient Algorithm with Momentum and Flexible Parameter Restart for Nonconvex Optimization}

\author{
Yi Zhou$^1$
\and
Zhe Wang$^2$\and
Kaiyi Ji$^{2}$\and
Yingbin Liang$^2$\And
Vahid Tarokh$^3$
\affiliations
$^1$University of Utah\\
$^2$The Ohio State University\\
$^3$Duke University
\emails
{yi.zhou}@utah.edu,
{wang.10982}@osu.edu,
{ji.367}@osu.edu
}

\begin{document}

\maketitle

\begin{abstract}
 Various types of parameter restart schemes have been proposed for proximal gradient algorithm with momentum to facilitate their convergence in convex optimization. However, under parameter restart, the convergence of proximal gradient algorithm with momentum remains obscure in nonconvex optimization. 
 In this paper, we propose a novel proximal gradient algorithm with momentum and parameter restart for solving nonconvex and nonsmooth problems. Our algorithm is designed to 1) allow for adopting flexible parameter restart schemes that cover many existing ones; 2) have a global sub-linear convergence rate in nonconvex and nonsmooth optimization; and 3) have guaranteed convergence to a critical point and have various types of asymptotic convergence rates depending on the parameterization of local geometry in nonconvex and nonsmooth optimization. Numerical experiments demonstrate the convergence and effectiveness of our proposed algorithm.   
\end{abstract}

\section{Introduction}
Training modern machine learning models in real applications typically involves highly nonconvex optimization, and some effective interesting examples include deep learning \cite{RELU}, nature language processing and computer vision, etc. To solve these nonconvex optimization problems, gradient-based algorithms \cite{Nesterov2014} are popular choices due to their simplicity, effectiveness as well as well-understood convergence guarantees. 

In practical training of machine learning models, momentum has been a successful and widely applied optimization trick that facilitates the convergence of gradient-based algorithms. Various types of momentum schemes have been developed, e.g., \cite{Nesterov2014,Beck2009,Tseng2010,Ghadimi2016b,Li:2015}, and have been shown to improve the order of convergence rates of gradient-based algorithms in solving convex and strongly convex optimization problems. In specific, gradient descent algorithms with momentum have been shown to achieve the complexity lower bound for convex optimization \cite{Nesterov2014,Beck2009} and have guaranteed convergence in nonconvex optimization \cite{Ghadimi2016b,Li:2015}. 

Despite the superior theoretical advantages of momentum acceleration schemes, they do not fully exploit the potential for acceleration. For example, the basic momentum scheme \cite{Nesterov2014,Beck2009} adopts a diminishing momentum coefficient for accelerating smooth convex optimization, and it does not provide much momentum acceleration after a large number of iterations. Also, for accelerating strongly convex optimization, the choice of momentum coefficient requires the knowledge of condition number of the Hessian matrix, which is typically unknown a priori. To resolve these issues and further facilitate the practical convergence of  gradient algorithms with momentum, various types of {\em parameter restart} techniques have been proposed, e.g., \cite{Donoghue2015,Fercoq2016,Fercoq2017,Giselsson2014,Kim2018,Liang2017,Lin2015,Liu2017,Renegar2018,Roulet2017}.  
In these works, it has been demonstrated that restarting algorithm parameters (i.e., variables and momentum coefficient) periodically can suppress the oscillations of the training loss induced by the extrapolation step and improve the practical convergence in {\em convex} optimization. In specific, parameter restart is typically triggered by certain occurrences that may slow down the convergence, such as function value divergence \cite{Donoghue2015,Renegar2018} and gradient mismatch \cite{Donoghue2015}, etc. Therefore, parameter restart can reduce the instability and oscillations caused by momentum. However, in {\em nonconvex} optimization, the applications of parameter restart to gradient algorithms with momentum require to deal with the following open issues. 

{\em (a)} While the convergence of gradient algorithms with momentum and parameter restart have been well explored in {\em convex} optimization, they are of lack of theoretical understandings in {\em nonconvex} optimization, which are important for modern machine learning purpose. {\em (b)} Previous works on gradient algorithms with momentum and restart for convex optimization are based on very specific restart schemes in order to have convergence guarantee, but practically the best restart scheme can be problem dependent. Therefore, it is much desired to design a momentum scheme that allows to adopt flexible parameter restart schemes with theoretical convergence guarantee.
{\em (c)} The existing gradient algorithms with momentum for nonconvex optimization have convergence guarantees at the cost of either introducing extra computation steps \cite{Li:2015,Li2017} or imposing restrictions on the objective function \cite{Ghadimi2016b}. It is important to explore whether parameter restart can help alleviate these costs or restrictions.

Considering all the issues above, we are motivated to design a gradient algorithm with momentum and parameter restart that {\em (a)} has convergence guarantee in nonconvex optimization, {\em (b)} allows to apply flexible restart schemes in practice and {\em (c)} avoids the existing weakness and restrictions in design of accelerated methods for nonconvex optimization. We summarize our contributions as follows.   

\subsection{Our Contributions}
We consider the problem of minimizing a smooth nonconvex function plus a (non)smooth regularizer. To solve such a class of problems, we propose APG-restart: a momentum-accelerated proximal gradient algorithm with parameter restart (see \Cref{alg: Acc-PGD}) and show that APG-restart satisfies the following properties.
\begin{itemize}[leftmargin=*,topsep=0pt]
	\item APG-restart allows for adopting any parameter restart scheme (hence covers many existing ones). In particular, it guarantees to make monotonic progress on function value between successive restart periods of iterations.
	\item The design of the proximal momentum component in APG-restart leverages the notion of generalized gradient mapping (see \cref{eq: grad_map}), which leads to convergence guarantee in nonconvex optimization. Also, APG-restart does not require extra computation steps compared to other accelerated algorithms for nonconvex optimization \cite{Li2017,Li:2015}, and removes the restriction of bounded domain on the regularizer function in existing works \cite{Ghadimi2016b}.   
	\item APG-restart achieves the stationary condition at a global sublinear convergence rate (see \Cref{lemma: Acc-PGD dynamic}).
	\item Under the Kurdyka-{\L}ojasiewicz (K{\L}) property of nonconvex functions (see \Cref{def: KL}), the variable sequence generated by APG-restart is guaranteed to converge to a critical point. Moreover, the asymptotic convergence rates of function value and variable sequences generated by APG-restart are fully characterized by the parameterization of the \KL property of the objective function. This work is the first study of gradient methods with momentum and parameter restart under the \KL property.
\end{itemize}

\subsection{Related Works}
\paragraph{Gradient algorithms with momentum and parameter restart:} Various types of parameter restart schemes have been proposed for accelerated gradient-based algorithms for convex optimization. Specifically, \cite{Donoghue2015} proposed to restart the accelerated gradient descent algorithm whenever certain function value-based criterion or gradient-based criterion is violated. These restart schemes were shown to achieve the optimal convergence rate without prior knowledge of the condition number of the function. \cite{Giselsson2014} further proposed an accelerated gradient algorithm with restart and established formal convergence rate analysis for smooth convex optimization. \cite{Lin2015} proposed a restart scheme that automatically estimates the strong convexity parameter and achieves a near-optimal iteration complexity. \cite{Fercoq2016,Fercoq2017} proposed a restart scheme for accelerated algorithms that achieves a linear convergence in convex optimization under the quadratic growth condition. \cite{Liu2017,Roulet2017} studied convergence rate of accelerated algorithms with restart in convex optimization under the error bound condition and the {\L}ojasiewicz condition, respectively. \cite{Renegar2018} proposed a restart scheme that is based on achieving a specified amount of decrease in function value. All these works studied accelerated gradient algorithms with restart in convex optimization, whereas this work focuses on nonconvex optimization. 


\paragraph{Nonconvex optimization under \KL property:} The Kurdyka-{\L}ojasiewicz property is a generalization of the {\L}ojasiewicz gradient inequality for smooth analytic functions to nonsmooth sub-analytic functions. Such a local property was then widely applied to study the asymptotic convergence behavior of various gradient-based algorithms in nonconvex optimization \cite{Attouch2009,Bolte2014,Zhou2016,Zhou_2017a}. The \KL property has also been applied to study convergence properties of accelerated gradient algorithms \cite{Li2017,Li:2015} and heavy-ball algorithms \cite{Ochs2018,Liang2016} in nonconvex optimization. Some other works exploited the \KL property to study the convergence of second-order algorithms in nonconvex optimization, e.g., \cite{Yi2018}.



\section{Preliminaries}
In this section, we introduce some definitions that are useful in our analysis later.
Consider a proper\footnote{An extended real-valued function $h$ is proper if its domain $\dom h := \{ x: h(x) < \infty \}$ is nonempty.} and lower-semicontinuous function $h:\mathbb{R}^d \to \mathbb{R}$ which is {\em not} necessarily smooth nor convex. We introduce the following generalized notion of derivative for the function $h$.
 
\begin{definition}(Subdifferential and critical point, \cite{vari_ana})\label{def:sub}
	The Frech\'et subdifferential $\widehat\partial h$ of function $h$ at $x\in \dom h$ is the set of $u\in \mathbb{R}^d$ defined as
	\begin{align*}
	\widehat\partial h(x) := \bigg\{u: \liminf_{z\neq x, z\to x} \frac{h(z) - h(x) - u^\intercal(z-x)}{\|z-x\|} \ge 0 \bigg\},
	\end{align*}
	and the limiting subdifferential $\partial h$ at $x\in\dom h$ is the graphical closure of $\widehat\partial h$ defined as:
	\begin{align*}
	\partial h(x) := \{ u: \exists x_k \to x, h(x_k) \to h(x), u_k \in \widehat{\partial} h(x_k) \to u \}.
	\end{align*}
	The set of critical points of $h$ is defined as $\mathbf{\crit}\!~h := \{ x: \zero\in\partial h(x) \}$. 
\end{definition}

Note that when the function $h$ is continuously differentiable, the limiting sub-differential $\partial h$ reduces to the usual notion of gradient $\nabla h$.
Next, we introduce the Kurdyka-{\L}ojasiewicz (K{\L}) property of a function $h$. Throughout, we define the distance between a point $x\in \mathbb{R}^d$ and a set $\Omega \subseteq \mathbb{R}^d$  as $\dist_\Omega(x) := \inf_{w\in \Omega} \|x - w\|$. 
\begin{definition}(\KL property, \cite{Bolte2014})\label{def: KL}
	A proper and lower-semicontinuous function $h$ is said to satisfy the \KL property if for every compact set $\Omega\subset \dom h$ on which $h$ takes a constant value $h_\Omega \in \mathbb{R}$, there exist $\varepsilon, \lambda >0$ such that for all $x \in \{z\in \mathbb{R}^d : \dist_\Omega(z)<\varepsilon, h_\Omega < h(z) <h_\Omega + \lambda\}$, the following inequality is satisfied
	\begin{align}\label{eq: KL}
	\varphi' \left(h(x) - h_\Omega\right) \dist_{\partial h(x)}(\zero) \ge 1,
	\end{align}
	where $\varphi'$ is the derivative of function $\varphi: [0,\lambda) \to \mathbb{R}_+$, which takes the form $\varphi(t) = \frac{c}{\theta} t^\theta$ for some $c>0, \theta\in (0,1]$.
\end{definition}

To elaborate, consider the case where $h$ is differentiable. Then, the \KL property in \cref{eq: KL} can be rewritten as 
\begin{align}
h(x) - h_\Omega \le C \|\nabla h(x)\|^{p} \label{eq: KLsimple}
\end{align}
for some constant $C>0$ and $p\in (1, +\infty)$. In fact, \Cref{eq: KLsimple} can be viewed as a generalization of the gradient dominance condition 
that corresponds to the special case of $p = 2$. A large class of functions have been shown to satisfy the \KL property, e.g., sub-analytic functions, logarithm and exponential functions, etc \cite{Bolte2007}. These function classes cover most of nonconvex objective functions encountered in practical applications, e.g., logistic loss, vector and matrix norms, rank, and polynomial functions, etc. Please refer to \cite[Section 5]{Bolte2014} and \cite[Section 4]{Attouch2010} for more example functions. 

To handle non-smooth objective functions, we introduce the following notion of proximal mapping.
\begin{definition}(Proximal mapping)\label{def:prox}
	For a proper and lower-semicontinuous function $h$, its proximal mapping at $x\in \mathbb{R}^d$ with parameter $\eta > 0$ is defined as:
	\begin{align}
	\prox{\eta h}(x) := \argmin_{z\in \mathbb{R}^d} \bigg\{h(z) + \frac{1}{2\eta}\|z - x\|^2\bigg\}.
	\end{align}
\end{definition}

\section{APG-restart for Nonsmooth \& Nonconvex Optimization}
In this section, we propose a novel momentum-accelerated proximal gradient with parameter restart (referred to as APG-restart) for solving nonsmooth and nonconvex problems. 

Consider the composite optimization problem of minimizing a smooth and nonconvex function $f:\mathbb{R}^d \to \mathbb{R}$ plus a possibly nonsmooth and convex function $g:\mathbb{R}^d \to \mathbb{R}$, which is written as 
\begin{align*}
	\min_{x\in \mathbb{R}^d} F(x):= f(x) + g(x). \tag{P}
\end{align*}

We adopt the following standard assumptions on the objective function $F$ in the problem (P).
\begin{assum}\label{assum: f+g}
	The objective function $F$ in the problem $\mathrm{(P)}$ satisfies:
	\begin{enumerate}[leftmargin=*]
		\item Function $F$ is bounded below, i.e.,  $F^*:=\inf_{x\in \mathbb{R}^d} F(x) > -\infty$;
		\item For any $\alpha \in \mathbb{R}$, the level set $\{x: F(x) \le \alpha\}$ is compact; 
		\item The gradient of $f$ is $L$-Lipschitz continuous and $g$ is lower-semicontinuous and convex.
	\end{enumerate}
\end{assum}

Under \Cref{assum: f+g}, we further introduce the following mapping for any $\eta>0$ and $x, u \in \mathbb{R}^d$:
\begin{align}
G_\eta(x,u) := \frac{1}{\eta} \big(x-\mathrm{prox}_{\eta g}(x-\eta u) \big). \label{eq: grad_map}
\end{align}
Such a mapping is well-defined and single-valued due to the convexity of $g$. Moreover,
the critical points of function $F$ (cf., \Cref{def:sub}) in the problem (P) can be alternatively characterized as $\mathbf{\crit}\!~F := \{ x: \zero\in G_\eta(x,\nabla f(x)) \}$. Therefore, $G_\eta(x,\nabla f(x))$ serves as a type of `gradient' at point $x$, and we refer to such a mapping as {\em gradient mapping} in the rest of the paper. In particular, the gradient mapping reduces to the usual notion of gradient when the nonsmooth part $g\equiv 0$.

\begin{algorithm}[H]
	\caption{APG-restart for nonconvex optimization}
	\label{alg: Acc-PGD}
	{\bf Input:} $K \in \mathbb{N}$, restart periods $q_0=0, \{q_{t}\}_{t\ge 1} \in \mathbb{N}$, stepsizes $\{\lambda_k\}_{k}, \{\beta_k\}_{k} >0.$
	
	{\bf Define:} $Q_t:= \sum_{\ell=0}^{t}q_\ell$.
	
	{\bf Initialize:} $x_{-1} \in \mathbb{R}^d$.
	
	\For{$k=0, 1, \ldots, K$}
	{
		Denote $t$ the largest integer such that $Q_t \le k$,\\
		Set: $\alpha_k = \frac{2}{k-Q_t+2}$,\\
		\If{$k= Q_t~$ for some $t\in \mathbb{N}$}
		{
			Reset: $x_{k} = y_k = x_{k-1},$
		}
		$z_{k} = (1-\alpha_{k+1})y_{k} + \alpha_{k+1} x_{k}$, \\
		$x_{k+1} =  x_k - \lambda_{k} G_{\lambda_k}(x_k, \nabla f(z_k))$, \\
		$y_{k+1} = z_{k} - \beta_{k} G_{\lambda_k}(x_k, \nabla f(z_k))$.	
	}	
	{\textbf{Output:} $x_K$.}
\end{algorithm}

To solve the nonsmooth and nonconvex problem (P), we propose the APG-restart algorithm that is presented in \Cref{alg: Acc-PGD}. APG-restart consists of new design of momentum schemes for updating the variables $x_k$ and $y_k$, the extrapolation step for updating the variable $z_k$ where $\alpha_{k+1}$ denotes the associated momentum coefficient, and the restart periods $\{q_t\}_t$. We next elaborate the two major ingredients of APG-restart: new momentum design and flexible restart scheduling with convergence guarantee.

\paragraph{New momentum design:} 
%
We adopt new momentum steps in APG-restart for updating the variables $x_k$ and $y_k$, which are different from those of the AG method in \cite{Ghadimi2016b} and we compare our update rules with theirs as follows.
\begin{equation}
\text{(APG-restart):}
\left\{
\begin{aligned}
\!x_{k+1} \!=\!  x_k \!-\! \lambda_{k} G_{\lambda_k}(x_k, \!\nabla f(z_k)), \\
\!y_{k+1} \!=\! z_{k} \!-\! \beta_{k} G_{\lambda_k}(x_k, \!\nabla f(z_k)).
\end{aligned}
\right.
\end{equation}
\begin{equation}
\text{(AG):}
\left\{
\begin{aligned}
\!x_{k+1} \!=\! \mathrm{prox}_{\lambda_k g}(x_k\!-\!\lambda_k \nabla f(z_k)) \\
\!y_{k+1} \!=\! \mathrm{prox}_{\lambda_k g}(z_k\!-\!\beta_k \nabla f(z_k))
\end{aligned}
\right.
\end{equation}
It can be seen from the above comparison that our APG-restart uses the same gradient mapping term $G_{\lambda_k}(x_k, \nabla f(z_k))$ to update both of the variables $x_k$ and $y_k$, while the AG algorithm in \cite{Ghadimi2016b} updates them using different proximal gradient terms. Consequently, our APG-restart is more computationally efficient as it requires to compute one gradient mapping per iteration while the AG algorithm needs to perform two proximal  updates.
On the other hand, the update rules of the AG algorithm guarantee convergence in nonconvex optimization only for functions of $g$ with bounded domain \cite{Ghadimi2016b}. Such a restriction rules out regularization functions with unbounded domain, which are commonly used in practical applications, e.g., $\ell_1, \ell_2$ regularization, elastic net, etc. In comparison, as we show in the analysis later, the update rules of APG-restart has guaranteed convergence in nonconvex optimization and does not require the regularizer $g$ to be domain-bounded.

\paragraph{Guarantee for any restart scheduling:} APG-restart retains the convergence guarantee with any restart scheduling. In specific, by specifying an arbitrary sequence of iteration periods $\{q_{t}\}_t \in \mathbb{N}$, APG-restart calls the restart operation at the end of each period (i.e., whenever $k=Q_t$ for some $t$). Upon restart, both $x_k$ and $y_k$ are reset to be the variable $x_{k-1}$ generated at the previous iteration, and the momentum coefficient $\alpha_k$ is reset to be 1. In the subsequent iterations, the momentum coefficient is diminished inversely proportionally to the number of iterations within the restart period. 


Since our APG-restart retains convergence guarantee for any restart periods $\{q_t \}_t$, it can implement any criterion that determines when to perform the parameter restart and have a convergence guarantee (see our analysis later). We list in \Cref{table: 1} some popular restart criteria from existing literature and compare their practical performance under our APG-restart framework in the experiment section later. We note that the restart criterion of the gradient mapping scheme implicitly depends on the gradient mapping, as $y_{k+1}-z_{k} \propto G_{\lambda_k}(x_k, \nabla f(z_k))$ from the update rule in \Cref{alg: Acc-PGD}.

\begin{table*}[ht]
	\setlength{\tabcolsep}{4pt}
	\center
	{\small 
	\begin{tabular}{ccccc}
		\toprule
		\begin{tabular}{@{}c@{}} Restart \\ scheme \end{tabular}  
		& \begin{tabular}{@{}c@{}} Fixed restart \\ \cite{Nesterov07}\end{tabular} 
		& \begin{tabular}{@{}c@{}} Function value  \\ \cite{Donoghue2015}\\\cite{Giselsson2014}\\\cite{Kim2018} \end{tabular} 
		& \begin{tabular}{@{}c@{}} Gradient mapping \\ \cite{Donoghue2015}\\\cite{Giselsson2014}\\\cite{Kim2018} \end{tabular}
		& \begin{tabular}{@{}c@{}} Non-monotonic \\ \cite{Giselsson2014} \end{tabular} \\ \midrule
		\begin{tabular}{@{}c@{}} Check \\ condition \end{tabular}
		& \begin{tabular}{@{}c@{}} $q_t\equiv q \in \mathbb{N}$ \\ for all $t$ \end{tabular}   
		& \begin{tabular}{@{}c@{}} restart whenever \\ $F(x_k)>F(x_{k-1})$ \end{tabular}
		& \begin{tabular}{@{}c@{}} restart whenever \\ $\inner{z_{k}\!-\!y_{k}}{y_{k+1}\!-\!z_{k}}$ \\ $\ge 0$  \end{tabular}
		& \begin{tabular}{@{}c@{}} restart whenever \\ $\inner{z_{k}\!-\!y_{k}}{y_{k+1}\!-\!\frac{z_{k}\!+\!x_k}{2}}$\\ $\ge 0$  \end{tabular} \\
		\bottomrule
	\end{tabular}
}
\caption{Restart conditions for different parameter restart schemes.}\label{table: 1}
\end{table*}

Performing parameter restart has appealing benefits. First, synchronizing the variables $x_k$ and $y_k$ periodically can suppress the deviation between them caused by the extrapolation step. This further helps to reduce the oscillation of the generated function value sequence. Furthermore, 
restarting the momentum coefficient $\alpha_k$ periodically injects more momentum into the algorithm dynamic, and therefore facilitates the practical convergence of the algorithm.

\section{Convergence Analysis of APG-restart}

In this section, we study the convergence properties of APG-restart in solving nonconvex and nonsmooth optimization problems. 
We first characterize the algorithm dynamic of APG-restart.
\begin{restatable}{lemma}{LemmaDynamicPGD}[Algorithm dynamic]\label{lemma: Acc-PGD dynamic}
	Let \Cref{assum: f+g} hold and apply \Cref{alg: Acc-PGD} to solve the problem (P). Set $\beta_k \equiv \frac{1}{8L}$ and $\lambda_k \in [\beta_k, (1+\alpha_{k+1})\beta_{k}]$. Then, the sequence $\{x_k \}_k$ generated by APG-restart satisfies: for all $t=1,2,...$
	\begin{align}
	F(x_{Q_t}) \le  &F(x_{Q_{t-1}}) - \frac{L}{4}\sum_{k=Q_{t-1}}^{Q_{t}-1} \|x_{k+1} - x_k\|^2, \label{eq: 9}\\
	\dist_{\partial F(x_{Q_t})}^2(\zero) &\le 162L^2 {\sum_{k=Q_{t-1}}^{Q_t-1}\|x_{k+1} - x_k\|^2}. \label{eq: 10}
	\end{align}
\end{restatable}

\Cref{lemma: Acc-PGD dynamic} characterizes the period-wise algorithm dynamic of APG-restart. In specific,  \cref{eq: 9} shows that the function value sequence generated by APG-restart is guaranteed to decrease between two adjacent restart checkpoint (i.e., $Q_{t-1}$ and $Q_t$), and the corresponding progress $F(x_{Q_{t-1}})- F(x_{Q_t})$ is bounded below by the square length of the iteration path between the restart checkpoints, i.e., $\sum_{k=Q_{t-1}}^{Q_t-1} \|x_{k+1} - x_k\|^2$. On the other hand, \cref{eq: 10} shows that the norm of the subdifferential at the $t$-th restart checkpoint is bounded by the square length of the same iteration path. In summary, the algorithm dynamic of APG-restart is different from that of traditional gradient-based algorithms in several aspects: First, the dynamic of APG-restart is characterized at the restart checkpoints, while the dynamic of gradient descent is characterized iteration-wise \cite{Attouch2009,Attouch2013}. As we elaborate later, such a property makes the convergence analysis of APG-restart more involved; Second, APG-restart makes monotonic progress on the function value between two adjacent restart checkpoints. In other accelerated gradient algorithms, such a monotonicity property is achieved by introducing a function value check step \cite{Li:2015} or an additional proximal gradient step \cite{Li2017}.    

Based on the algorithm dynamic in \Cref{lemma: Acc-PGD dynamic}, we obtain the following global convergence rate of APG-restart for nonconvex and nonsmooth optimization. Throughout the paper, we denote $f(n) = \Theta(g(n))$ if and only if for some $0<c_1<c_2$, $c_1 g(n) \le f(n) \le c_2 g(n)$ for all $n\ge n_0$.  

\begin{restatable}{thm}{TheoremGlobal}[Global convergence rate]\label{thm: global}
	Under the same conditions as those of \Cref{lemma: Acc-PGD dynamic}, the sequence $\{z_k \}_k$ generated by APG-restart satisfies: for all $K=1,2,...$
\begin{align*}
\min_{0\le k\le K-1}\|G_{\lambda_k}(z_k, \nabla f(z_k))\|^2 \le \Theta \Big({\frac{L\big(F(x_{0}) - F^*\big)}{K}}\Big).
\end{align*}
\end{restatable}
\Cref{thm: global} establishes the global convergence rate of APG-restart in terms of the gradient mapping, which we recall characterizes the critical point of the nonconvex objective function $F$. In particular, the order of the above global convergence rate matches that of other accelerated gradient algorithms \cite{Ghadimi2016b} for nonconvex optimization, and APG-restart further benefits from the flexible parameter restart scheme that provides extra acceleration in practice (as we demonstrate via experiments later).   

\Cref{thm: global} does not fully capture the entire convergence property of APG-restart. To elaborate, convergence of the gradient mapping in \Cref{thm: global} does not necessarily guarantee the convergence of the {\em variable sequence} generated by APG-restart. On the other hand, the convergence rate estimate is based on the global Lipschitz condition of the objective function, which may not capture the local geometry of the function around critical points and therefore leads to a coarse convergence rate estimate in the asymptotic regime. To further explore stronger convergence results of APG-restart, we next exploit the ubiquitous Kurdyka-{\L}ojasiewicz (K{\L}) property (cf., \Cref{def: KL}) of nonconvex functions. We make the following assumption.

\begin{assum}\label{assum: KL}
	The objective function $F$ in the problem $\mathrm{(P)}$ satisfies the \KL property.
\end{assum}

Based on the algorithm dynamic in \Cref{lemma: Acc-PGD dynamic} and further leveraging the \KL property of the objective function, we obtain the following convergence result of APG-restart in nonconvex optimization.

\begin{restatable}{thm}{TheoremVariable}[Variable convergence]\label{thm: Acc-GD variable}
	Let Assumptions \ref{assum: f+g} and \ref{assum: KL} hold and apply \Cref{alg: Acc-PGD} to solve the problem (P). Set $\beta_k \equiv \frac{1}{8L}$ and $\lambda_k \in [\beta_k, (1+\alpha_{k+1})\beta_{k}]$. Define the length of iteration path of the $t$-th restart period as
	$	L_t:= \sqrt{\sum_{k=Q_t}^{Q_{t+1}-1} \|x_{k+1}-x_k\|^2}.$
	Then, the sequence $\{L_t \}_t$ generated by APG-restart satisfies: for all periods of iterations $t=1,2,...$
	\begin{align}
	\sum_{t=0}^\infty L_t < +\infty. \label{eq: finite_len}
	\end{align}
	Consequently, the variable sequences $\{x_k \}_k, \{y_k \}_k, \{z_k \}_k$ generated by APG-restart converge to the same critical point of the problem (P), i.e., 
	\begin{align}
		x_k,y_k,z_k \overset{k}{\to} x^* \in \mathbf{\crit}\!~F.
	\end{align}
\end{restatable}

\Cref{thm: Acc-GD variable} establishes the formal convergence of APG-restart in nonconvex optimization. We note that such a convergence guarantee holds for any parameter restart schemes, therefore demonstrating the flexibility and generality of our algorithm. Also, unlike other accelerated gradient-type of algorithms that guarantee only convergence of function value \cite{Li:2015,Li2017}, our APG-restart is guaranteed to generate convergent variable sequences to a critical point in nonconvex optimization. 

To highlight the proof technique, we first exploit the dynamic of APG-restart in \Cref{lemma: Acc-PGD dynamic} to characterize the limit points of the sequences $\{x_{Q_t} \}_t, \{F(x_{Q_t}) \}_t$ that are indexed by the restart checkpoints. Then, we further show that the entire sequences  $\{x_{k} \}_k, \{F(x_{k}) \}_k$ share the same limiting properties, which in turn guarantee the sequences to enter a local parameter region of the objective function where the \KL property can be exploited. Taking advantage of the \KL property, we are able to show that the length of the optimization path is finite as iteration $k\to \infty$. Consequently, the generated variable sequences can be shown to  converge to a certain critical point of the Problem (P). 

Besides the variable convergence guarantee under the \KL property, we also obtain various types of convergence rate estimates of APG-restart depending on the specific parameterization of the local \KL property of the objective function. We obtain the following results. 

\begin{restatable}{thm}{TheoremRates}[Convergence rate of function value]\label{thm: Acc-GD rates}
	Let Assumptions \ref{assum: f+g} and \ref{assum: KL} hold and apply \Cref{alg: Acc-PGD} to solve the problem (P). Set $\beta_k \equiv \frac{1}{8L}$ and $\lambda_k \in [\beta_k, (1+\alpha_{k+1})\beta_{k}]$. Suppose the algorithm generates a sequence $\{x_k \}_k$ that converges to a certain critical point $x^*$ where the K{\L} property holds with parameter $\theta\in (0,1]$. Then, there exists a sufficiently large $t_0\in \mathbb{N}$ such that for all $t\ge t_0$,
	\begin{enumerate}[leftmargin=*]
		\item If $\theta = 1$, then $F(x_{Q_t}) \downarrow F(x^*)$ within finite number of periods of iterations; 
		\item If $\theta \in [\frac{1}{2}, 1)$, then $F(x_{Q_t}) \downarrow F(x^*)$ linearly as 
		$F(x_{Q_t}) - F(x^*) \le \exp \big(-\Theta(t-t_0)\big);$
		\item If $\theta \in (0, \frac{1}{2})$, then $F(x_{Q_t}) \downarrow F(x^*)$ sub-linearly as 
		$F(x_{Q_t}) - F(x^*) \le \Theta \Big((t-t_0)^{-\frac{1}{1-2\theta}}\Big).$
	\end{enumerate}
\end{restatable}

\begin{restatable}{thm}{TheoremRatesVariable}[Convergence rate of variable]\label{thm: Acc-GD var_rates}
	Under the same conditions as those of \Cref{thm: Acc-GD rates}, suppose APG-restart generates a sequence $\{x_k \}_k$ that converges to a certain critical point $x^*$ where the K{\L} property holds with parameter $\theta\in (0,1]$. Then, there exists a sufficiently large $t_0\in \mathbb{N}$ such that for all $t\ge t_0$,
	\begin{enumerate}[leftmargin=*]
		\item If $\theta = 1$, then $x_{Q_t} \overset{t}{\to} x^*$ within finite number of periods of iterations; 
		\item If $\theta \in [\frac{1}{2}, 1)$, then $x_{Q_t} \overset{t}{\to} x^*$ linearly as 
		$\|x_{Q_t} - x^*\| \le \exp \big(-\Theta(t-t_0)\big);$
		\item If $\theta \in (0, \frac{1}{2})$, then $x_{Q_t} \overset{t}{\to} x^*$ sub-linearly as $\|x_{Q_t} - x^*\| \le\Theta \Big((t-t_0)^{-\frac{\theta}{1-2\theta}}\Big)$.
	\end{enumerate}
\end{restatable}

\Cref{thm: Acc-GD rates} and \Cref{thm: Acc-GD var_rates} establish the asymptotic convergence rate results for the function value sequence and variable sequence generated by APG-restart, respectively. Intuitively, after a sufficiently large number of training iterations, APG-restart enters a local neighborhood of a certain critical point. In such a case, the global convergence rate characterized in \Cref{thm: global} can be a coarse estimate because it exploits only the global Lipschitz property of the function. On the contrary, the local \KL property characterizes the function geometry in a more accurate way and leads to the above tighter convergence rate estimates. In particular, the \KL parameter $\theta$ captures the `sharpness' of the local geometry of the function, i.e., a larger $\theta$ induces a faster convergence rate.

\section{Experiments}\label{sec:exp}
In this section, we implement the APG-restart algorithm with different restart schemes listed in \Cref{table: 1} to corroborate our theory that APG-restart has guaranteed convergence with any restart scheme. In specific, for the fixed restart scheme we set the restart period to be $q=10,30,50$, respectively. 

We first solve two smooth nonconvex problems, i.e., the logistic regression problem with a nonconvex regularizer (i.e., $g(x):=\alpha \sum_{i=1}^{d} \frac{x_i^2}{1+x_i^2}$) and the robust linear regression problem. For the logistic regression problem, we adopt the cross-entropy loss and set $\alpha=0.01$, and for the robust linear regression problem, we adopt the robust nonconvex loss $\ell(s):= \log(\frac{s^2}{2}+1)$. We test both problems on two LIBSVM datasets: a9a and w8a \cite{Chang_2011}. We use stepsizes $\beta_k = 1, \lambda_{k}=(1+\alpha_{k+1})\beta_{k}$ for the APG-restart as suggested by our theorems. We note that in these experiments, we plot the loss gap versus number of iterations for all algorithms. The comparison of running time is similar as all the algorithms require the same computation per iteration.

\begin{figure} 
	\begin{subfigure}{0.48\linewidth}
		\includegraphics[width=\linewidth]{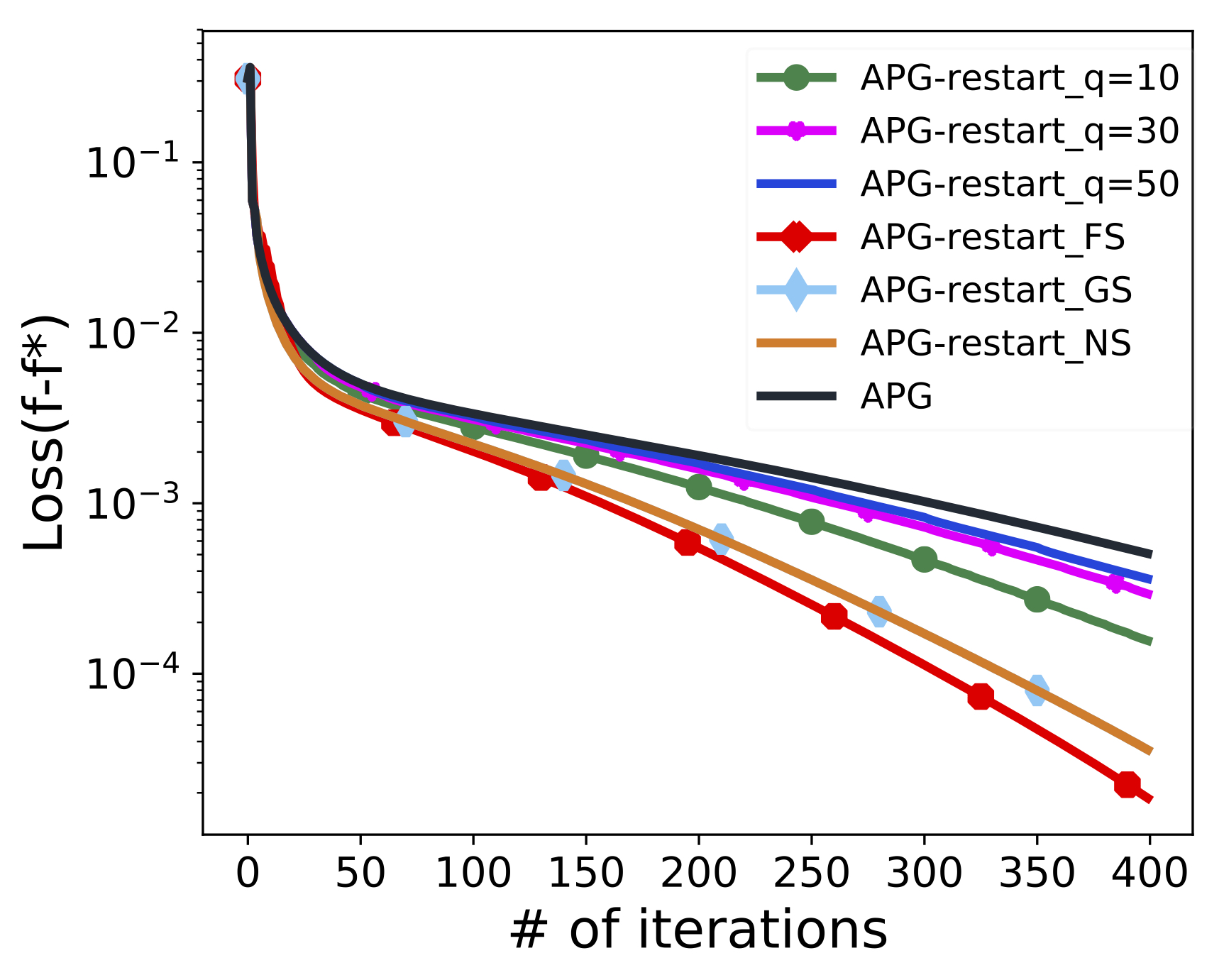}
		\caption{logistic regression \\ a9a}
	\end{subfigure}
	\begin{subfigure}{0.48\linewidth}
		\includegraphics[width=\linewidth]{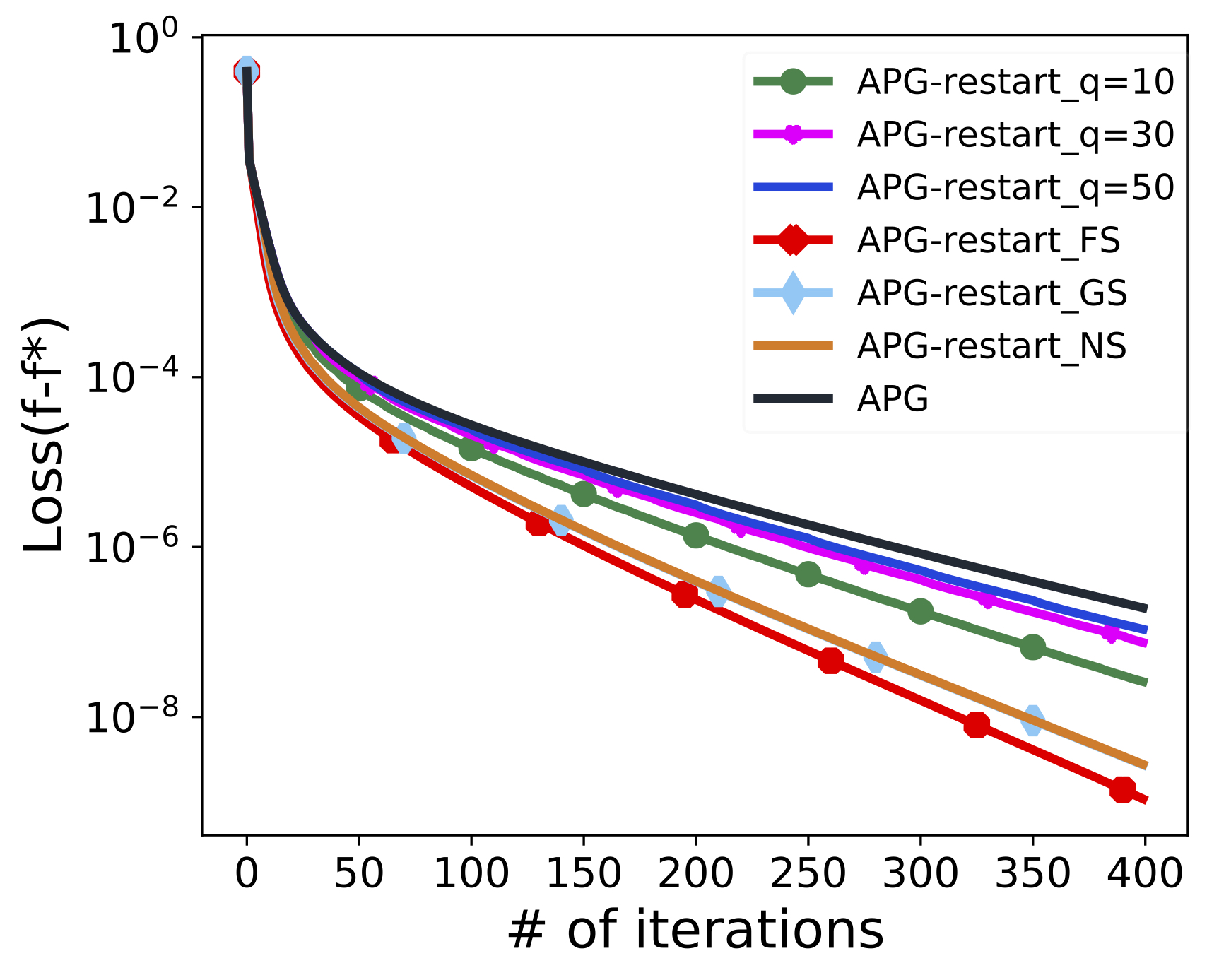}
		\caption{logistic regression \\ w8a}
	\end{subfigure}%
	\\
\begin{subfigure}{0.49\linewidth}
	\includegraphics[width=\linewidth]{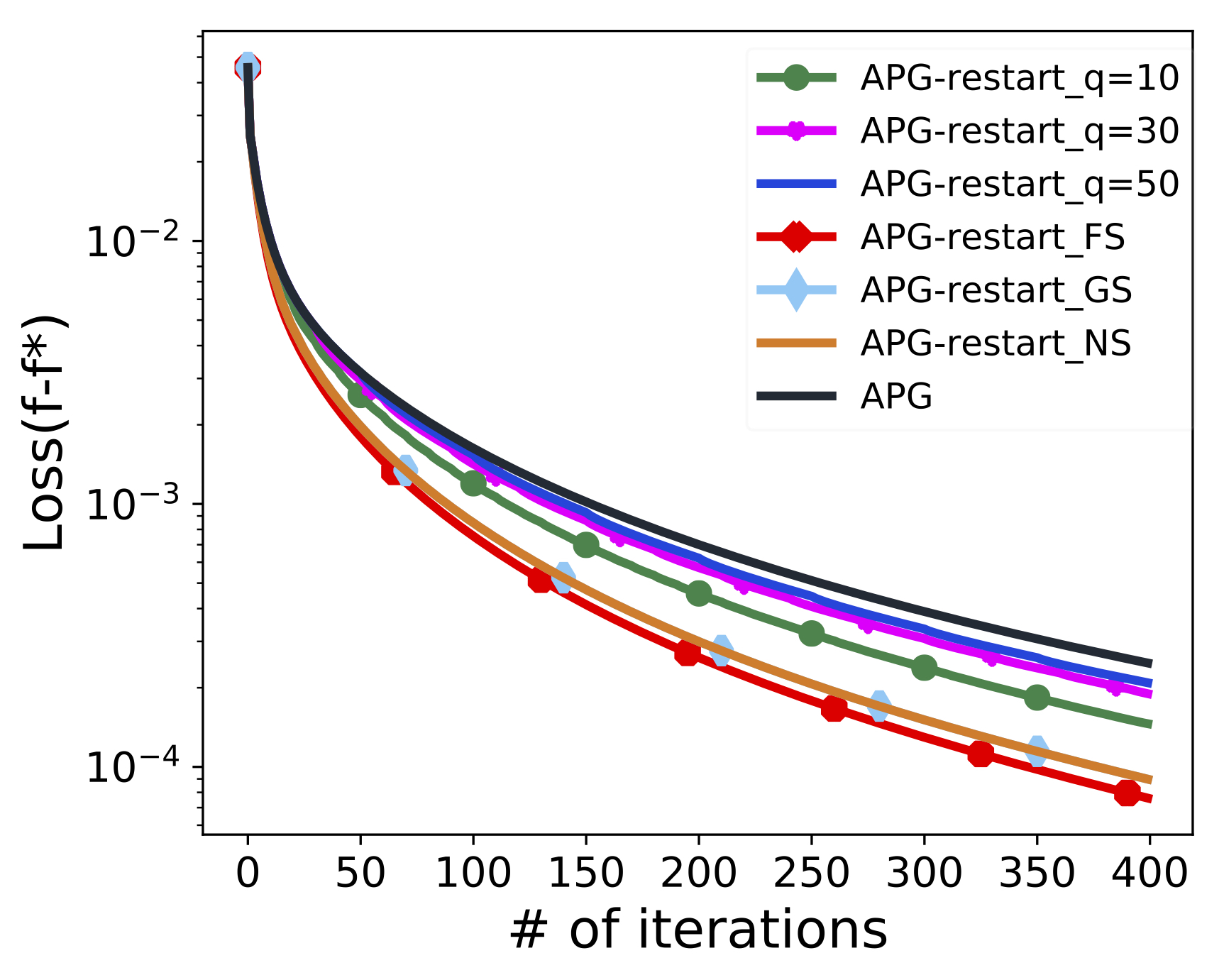}
	\caption{robust regression \\ a9a}
\end{subfigure}
\begin{subfigure}{0.49\linewidth}
	\includegraphics[width=\linewidth]{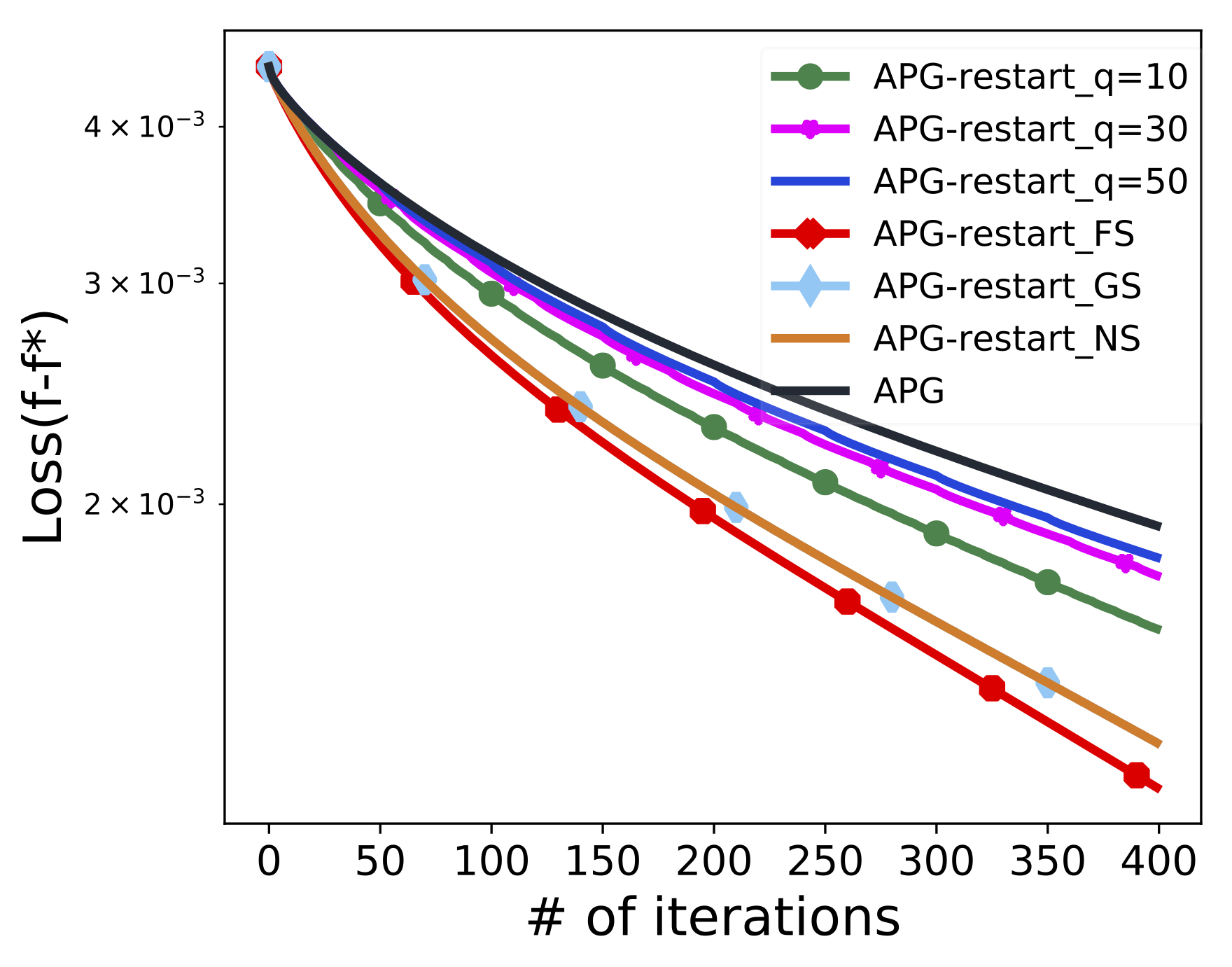}
	\caption{robust regression \\ w8a}
\end{subfigure}%
	\caption{Comparison of different restart schemes in smooth nonconvex optimization.}   \label{Experment_1}
\end{figure}

\Cref{Experment_1} shows the experiment results of APG-restart with fixed scheme (constant $q$), function value scheme (FS), gradient mapping scheme (GS) and non-monotone scheme (NS). It can be seen that APG-restart under the function scheme performs the best among all restart schemes. In fact, the function scheme restarts the APG algorithm the most often in these experiments. The gradient mapping scheme and the non-monotone scheme have very similar performance, and both of them perform slightly worse than the function scheme. Moreover, the fixed restart schemes have the worst performance. In particular, the performance of fixed scheme gets better as the restart period $q$ decreases (i.e., more restarts take place).    

\begin{figure} 
	\centering 
	\begin{subfigure}{0.48\linewidth}
		\includegraphics[width=\linewidth]{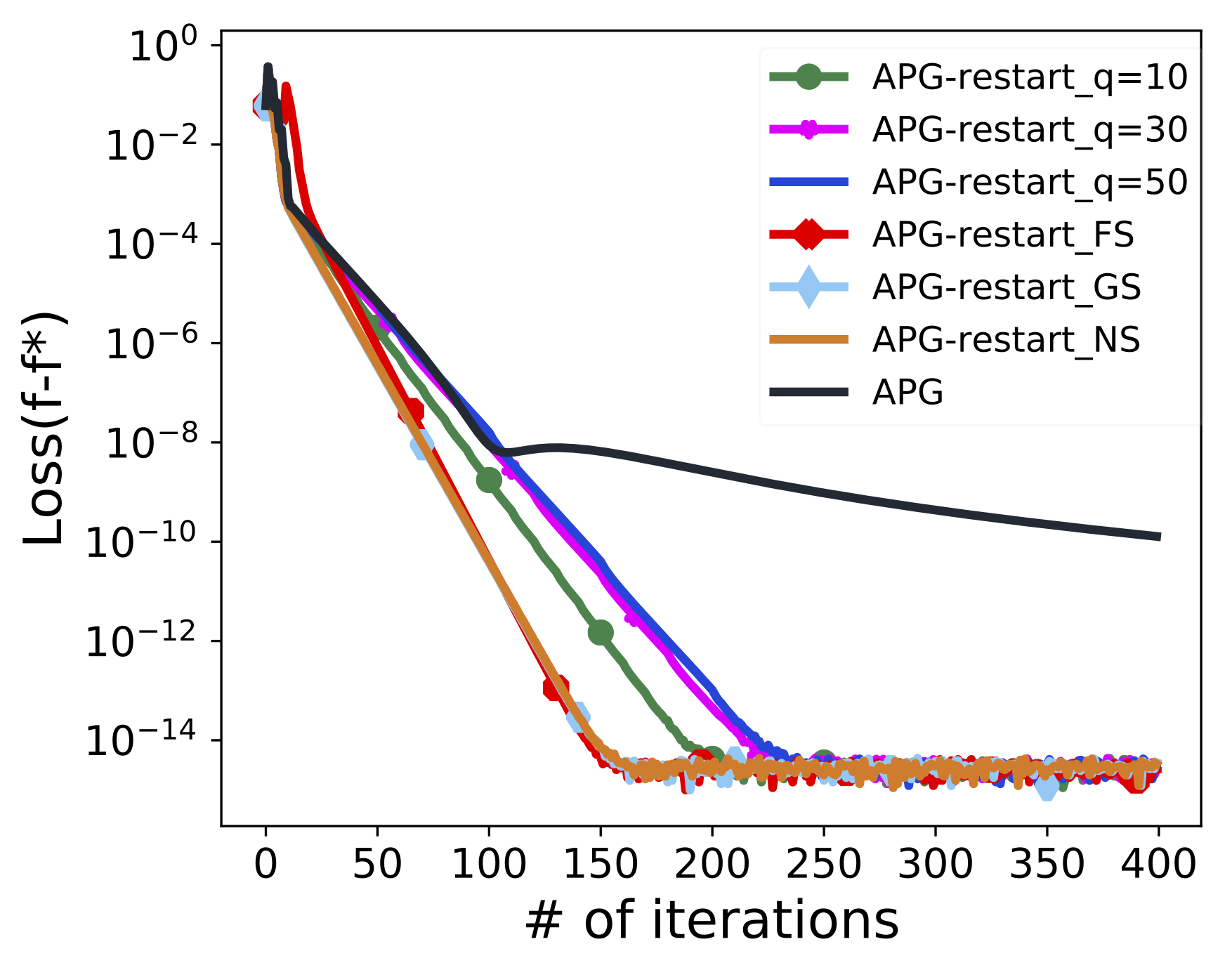}
		\caption{logistic regression \\a9a}
	\end{subfigure}
	\begin{subfigure}{0.48\linewidth}
		\includegraphics[width=\linewidth]{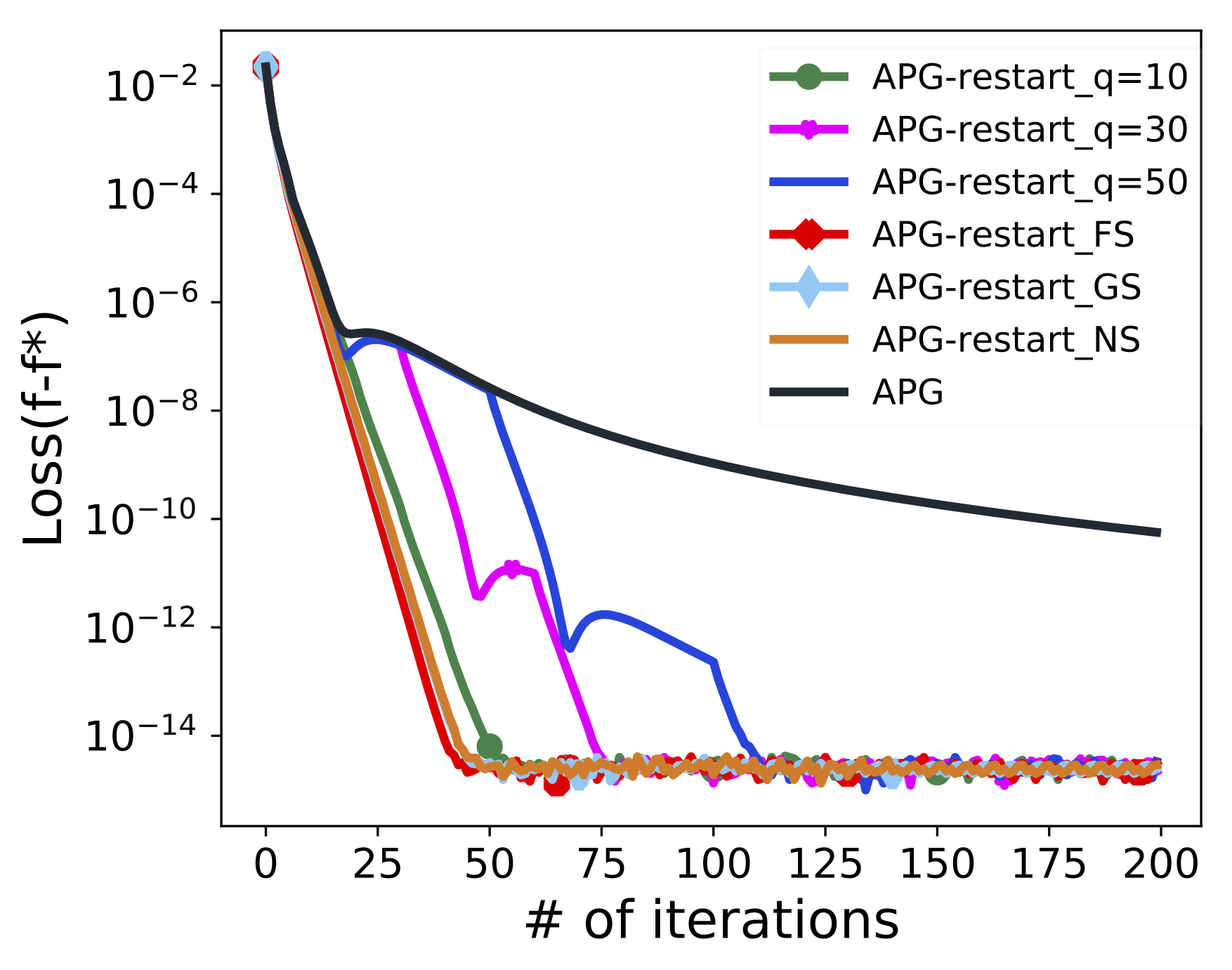}
		\caption{logistic regression \\ w8a}
	\end{subfigure}%
	\\
\begin{subfigure}{0.48\linewidth}
	\includegraphics[width=\linewidth]{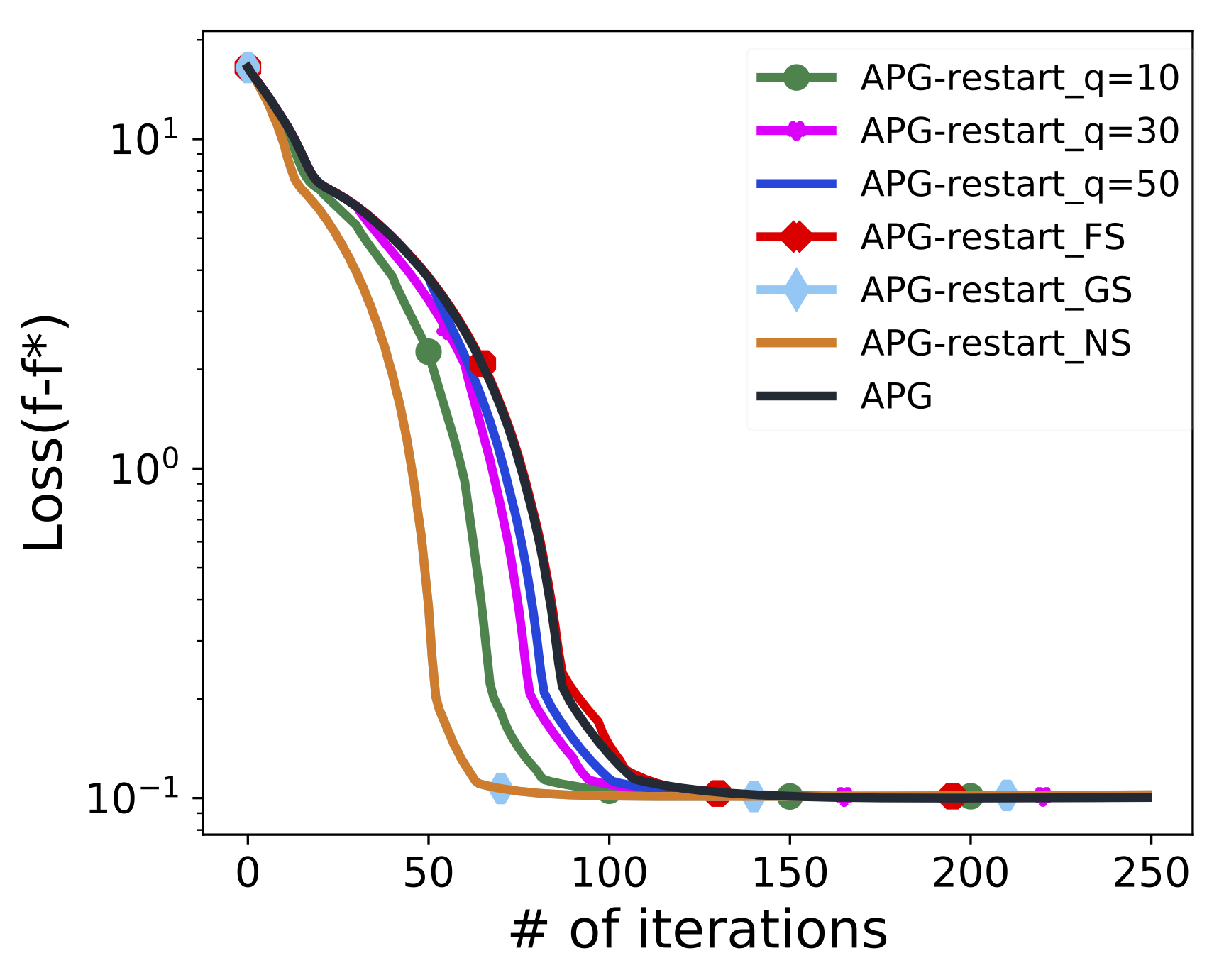}
	\caption{robust regression \\ a9a}
\end{subfigure}
\begin{subfigure}{0.48\linewidth}
	\includegraphics[width=\linewidth]{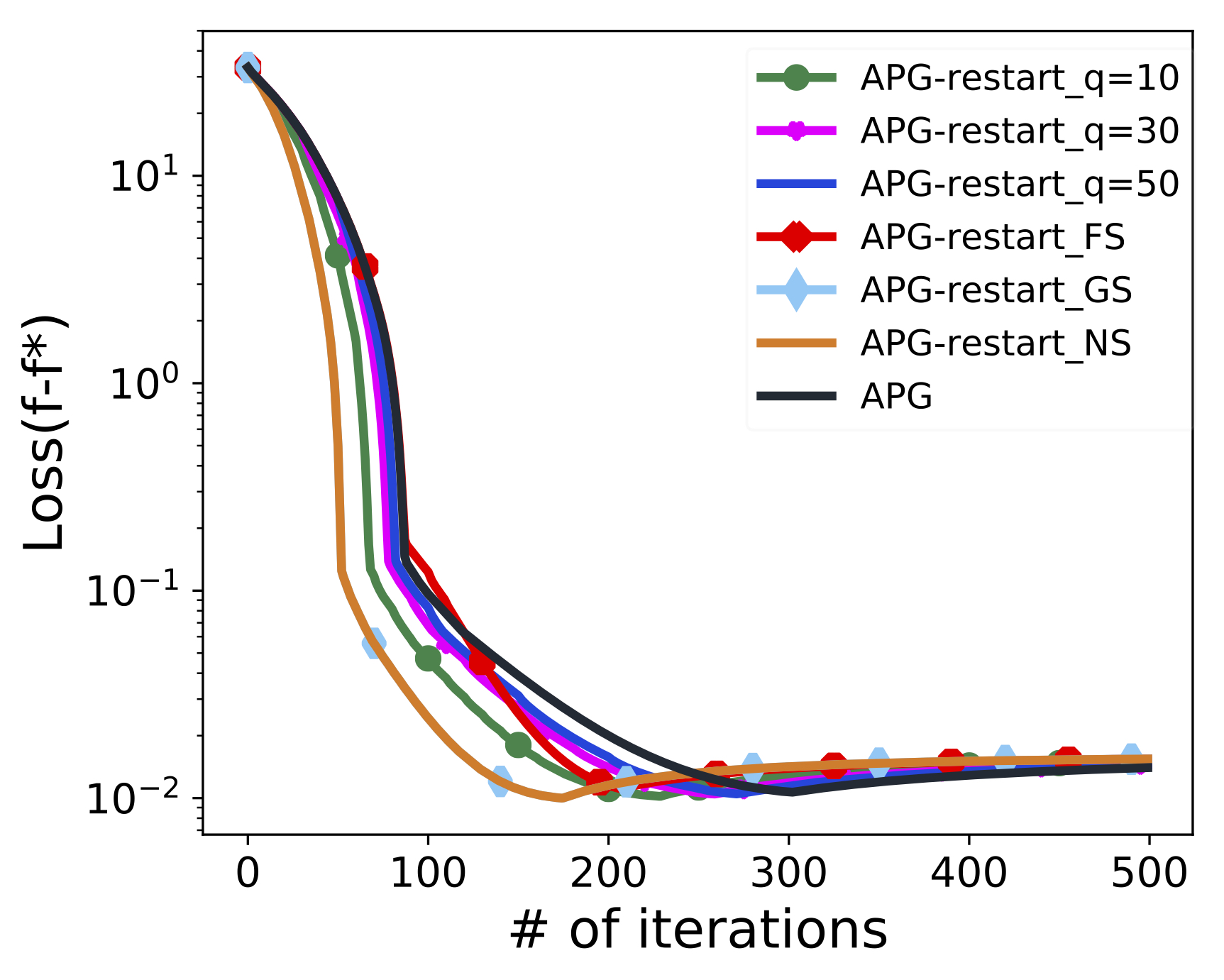}
	\caption{robust regression \\ w8a}
\end{subfigure}%
	\caption{Comparison of different restart schemes in {\em nonsmooth} nonconvex optimization.}   \label{Experment_3}
\end{figure}

Next, we further add a nonsmooth $\ell_1$ norm regularizer to the objective functions of all the problems mentioned above, and apply APG-restart with different restart schemes to solve them. The results are shown in \Cref{Experment_3}. One can see that for the nonsmooth logistic regression, all the non-fixed restart schemes have comparable performances and they perform better than the fixed restart schemes. For the nonsmooth robust linear regression, both the gradient mapping scheme and the non-monotone scheme outperform the other schemes. In this experiment, the function scheme has a degraded performance that is comparable to the fixed restart schemes. This is possibly due to the highly nonconvexity of the loss landscape.

\section{Conclusion}
In this paper, we propose a novel accelerated proximal gradient algorithm with parameter restart for nonconvex optimization. Our proposed APG-restart allows for adopting any parameter restart schemes and have guaranteed convergence. We establish both the global convergence rate and various types of asymptotic convergence rates of the algorithm, and we demonstrate the effectiveness of the proposed algorithm via numerical experiments. We expect that such a parameter restart algorithm framework can inspire new design of optimization algorithms with faster convergence for solving nonconvex machine learning problems.

{\small
	\section*{Acknowledgment}
	The work of Z. Wang, K. Ji and Y. Liang was supported in part by the U.S. National Science Foundation under the grants CCF-1761506, CCF-1909291 and CCF-1900145.
	
	\bibliographystyle{named}
	\bibliography{./ref}
}

\end{document}